\PassOptionsToPackage{colorlinks=true, allcolors=blue}{hyperref}
\documentclass[10pt]{article}
\usepackage{authblk}
\usepackage{orcidlink}
\usepackage{graphicx}
\usepackage{subcaption}

\usepackage{amsmath, amsthm, amssymb}
\usepackage{marvosym}
\usepackage[a4paper]{geometry}
\usepackage{comment}

\usepackage[colorlinks=true, allcolors=blue]{hyperref}

\title{\textbf{Retrieving biparameter persistence modules from monoparameter ones: a characterization of hook-decomposable persistence modules}}
\author[1,2]{Isabella Mastroianni\,\orcidlink{0009-0002-9866-3648}}
\author[1]{Marco Guerra\,\orcidlink{0000-0003-0033-3748}}
\author[1]{Ulderico Fugacci\,\orcidlink{0000-0003-3062-997X}}
\author[2]{Emanuela De Negri\,\orcidlink{0000-0001-5556-440X}}

\affil[1]{\textit{Institute of Applied Mathematics and Information Technologies “Enrico Magenes”, National Research Council, Genova, Italy}}
\affil[2]{\textit{Department of Mathematics, University of Genova, Italy}}
\date{June 2025}

\newtheorem{proposition}{Proposition}
\newtheorem{oss}{Remark}

\newtheorem{definition}{Definition}
\newtheorem{example}{Example}

\begin{document}

\maketitle

\begin{abstract}
    Motivated by the need to relate the biparameter persistence module induced by a pair of scalar functions with the monoparameter persistence modules induced by each function separately, we introduce a construction \, ---\, the \(\gamma\)-product \, ---\, that defines a kind of product between two monoparameter persistence modules.

    While originally conceived to serve this comparative purpose, our construction unexpectedly reveals a deeper structural property: it also characterizes a class of biparameter modules known as hook-decomposable modules.
\end{abstract}

Due to space limitations, proofs are omitted and will appear in a forthcoming extended version of the work.

Let \( f, g \colon X \to \mathbb{N} \) be continuous functions defined on a topological space \( X \). Fix \( k \in \mathbb{N} \) and a field \(\mathbb{F}\), and consider the two monoparameter persistence modules \( M_f \) and \( M_g \), as well as the biparameter persistence module \( M_{(f,g)} \), obtained by applying \( H_k(\cdot\thinspace; \mathbb{F}) \), to the filtrations induced by \( f \), \( g \), and the pair \( (f, g) \), respectively \cite{edelsbrunner2002topological, zomorodian2004computing, carlsson2007theory}.

We work under the assumption that all persistence modules involved are of finite type. Moreover, we identify them with their corresponding finitely generated graded modules over polynomial rings \cite{zomorodian2004computing, oudot2017persistence}. In particular, \( M_f \) is an \( \mathbb{N} \)-graded \( \mathbb{F}[x] \)-module, \( M_g \) is an \( \mathbb{N} \)-graded \( \mathbb{F}[y] \)-module, and \( M_{(f,g)} \) is an \( \mathbb{N}^2 \)-graded \( \mathbb{F}[x,y] \)-module.

By \(PD_k(f)\), we denote the \( k \)-th persistence diagram of \( f \). We assume \(PD_k(f)\) contains infinitely many copies of the points in the diagonal \(\Delta\) \cite{edelsbrunner2002topological}. The same holds for \( g \).

In order to work with sets of generators for \(M_f\) and \(M_g\) which are both of the same finite cardinality, we define the following sets depending on a bijection \( \gamma\colon PD_k(f) \to PD_k(g) \). 
We define \(S_\gamma\) as:
\[
PD_k(f) \setminus \{(b,b) \in PD_k(f) \mid \gamma(b,b) \in \Delta\}.
\]

By construction of \(S_\gamma\), a presentation of \( M_f \) is:
\[
\left\langle m_{(b,d)} \;\middle|\; (b,d) \in S_\gamma, \; \deg(m_{(b,d)}) = b,\; m_{(b,d)} x^{d-b} = 0 \right\rangle.
\]
Similarly, by considering \( S_{\gamma^{-1}} \), we obtain an analogous presentation of \( M_g \) with the same number of generators \( m'_{(b',d')} \) and relations \( m'_{(b',d')} y^{d'-b'} = 0 \), where we denote \(\gamma(b,d) \) by \( (b', d')\).

We are now ready to state the core definition of this work.

\begin{definition}[$\gamma$-product of \( M_f \) and \( M_g \)]
Let \( M_f \) and \( M_g \) be as above. We define the \emph{\( \gamma \)-product} of \( M_f \) and \( M_g \) as the $\mathbb{N}^2$-graded \( \mathbb{F}[x,y] \)-module \( M_f \MVAt_\gamma M_g \) with the following structure:
\begin{itemize}
  \item generators are the pairs \( (m_{(b,d)},\, m'_{(b',d')}) \), each of bidegree \( (b, b') \);
  \item relations are given by:
    \[
    (m_{(b,d)},\, m'_{(b',d')}) \cdot x^{d - b} y^{d' - b'} = 0;
    \]
  \item addition is defined component-wise: \[ (m, m') + (n, n') := (m+n,\, m'+n') ;
  \]
  \item \( \mathbb{F}[x,y] \)-action and scalar multiplication are concisely defined by: 
    \[
    (m,m') \cdot s(x,y):=\big(m\cdot s(x,1),\, m'\cdot s(1,y)\big).
    \]
\end{itemize}
\end{definition}

\begin{proposition}
\( M_f \MVAt_{\gamma} M_g \) is a finitely generated \(\mathbb{N}^2 \)-graded \( \mathbb{F}[x,y] \)-module.
\end{proposition}

    The proof is a straightforward computation, which we omit here, that relies on the graded \( \mathbb{F}[x] \)- and \( \mathbb{F}[y] \)-module structures of \( M_f \) and \( M_g \). Before focusing on an explicit example, let us examine some basic properties of this construction.

\begin{oss}
Although the construction of \( M_f \MVAt_{\gamma} M_g \) might resemble a tensor product, it is important to notice that it does not define a tensor structure in the categorical sense. In particular, this construction is not bilinear; instead, it should be understood as a matching-based combination of two persistence modules, whose structure depends explicitly on the chosen \( \gamma \).
\end{oss}


We now introduce a bijection that is particularly relevant to our framework.

\begin{definition}
We define \( \overline{\gamma} \) as a bijection minimizing the interleaving distance:
\[
d_I(M_{(f,g)}, M_f \MVAt_{\overline{\gamma}} M_g) = \inf_{\gamma} d_I(M_{(f,g)}, M_f \MVAt_{\gamma} M_g).
\]
\end{definition}

Notice that in general \( \overline{\gamma} \) does not realize the bottleneck distance between the persistence diagrams of \( f \) and \( g \), as illustrated in the following example.

\begin{example}\label{ex:gamma bar}
Consider the two functions \( f \) and \( g \) defined simplex-wise on the simplicial complex \( X \) as shown in Figure~\ref{fig:X and f,g}.

\begin{figure}[!htbp]
    \centering
    \includegraphics[width=0.4\textwidth]{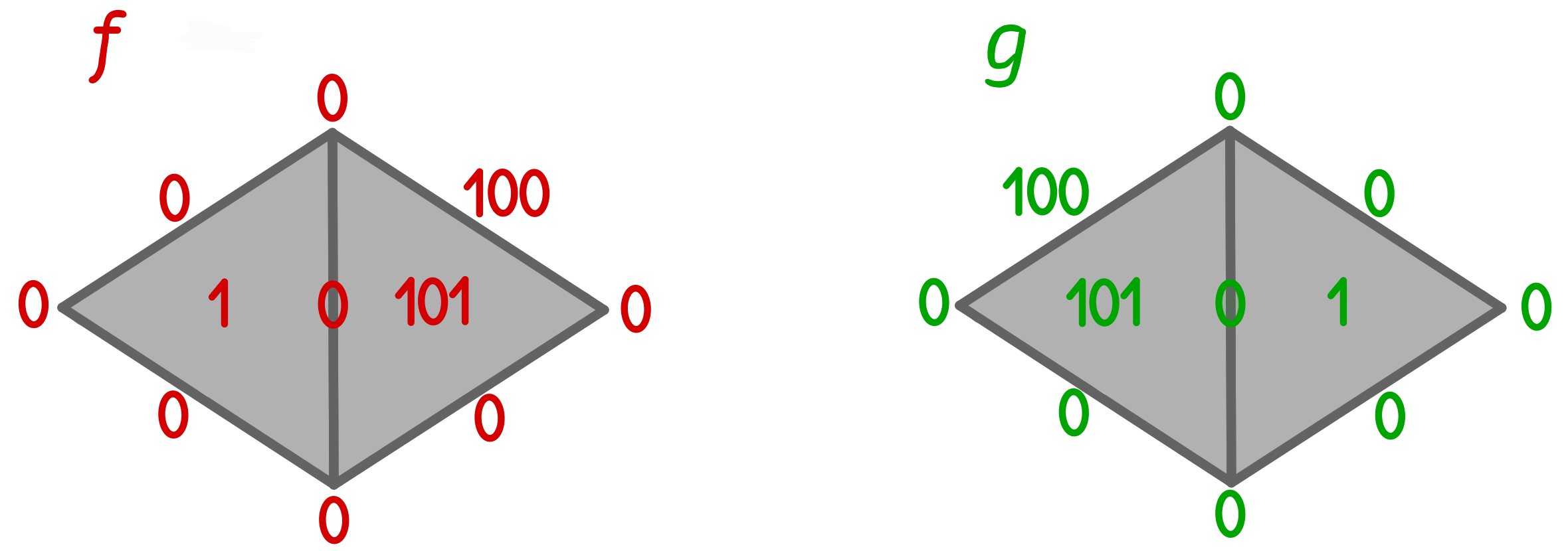}
    \caption{A simplicial complex \( X \) and two functions \( f, g \) defined simplex-wise on it.}
    \label{fig:X and f,g}
\end{figure}

The first homology of the filtrations induced by \( f \) and \( g \) yields two identical persistence diagrams, each containing the pairs \( (0,1) \) and \( (100,101) \). Hence, the identity bijection \( \gamma_{\mathrm{bott}} \) realizes the bottleneck distance. The corresponding product \( M_f \MVAt_{\gamma_{\mathrm{bott}}} M_g \) has the support shown in Figure~\ref{fig:subfig-supp Mf bott Mg}.

However, the support of \( M_{(f,g)} \), shown in Figure~\ref{fig:subfig-supp Mfg},  differs.

\begin{figure}[!htbp]
\centering
\begin{subfigure}[b]{0.43\textwidth}
    \centering
    \includegraphics[width=\textwidth]{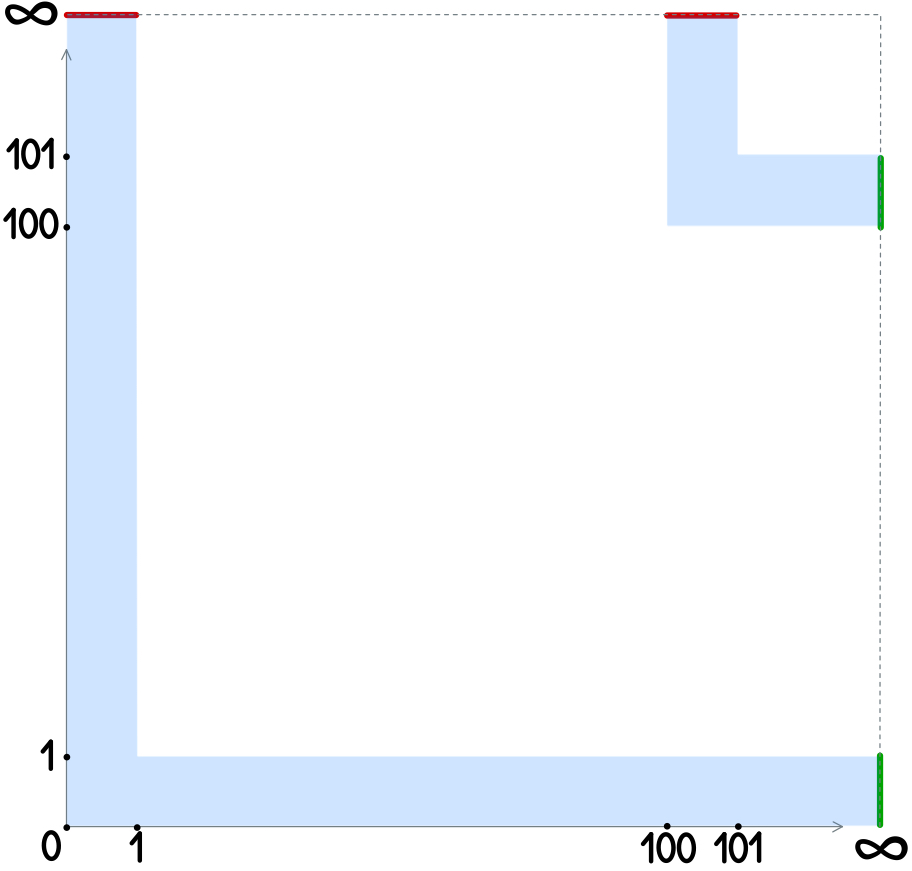}
    \caption{\footnotesize Support of \( M_f \MVAt_{\gamma_{\mathrm{bott}}} M_g \).}
    \label{fig:subfig-supp Mf bott Mg}
\end{subfigure}
\hfill
\begin{subfigure}[b]{0.43\textwidth}
    \centering
    \includegraphics[width=\textwidth]{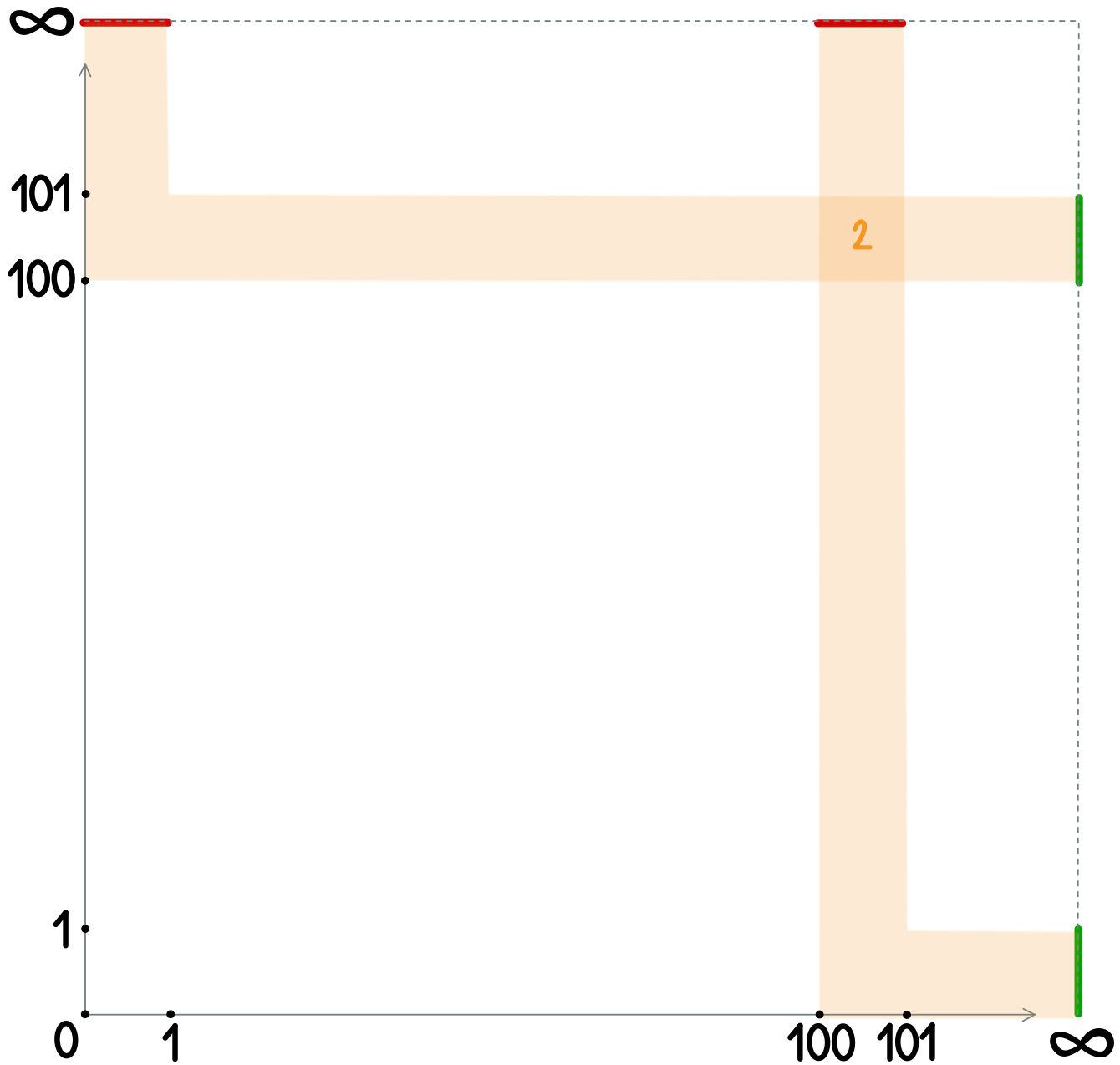}
    \caption{\footnotesize Support of \( M_{(f,g)} \).}
    \label{fig:subfig-supp Mfg}
\end{subfigure}
\caption{Comparison of the supports of (a) \( M_f \MVAt_{\gamma_{\mathrm{bott}}} M_g \) and (b) \( M_{(f,g)} \).}
\label{fig:supp-comparison}
\end{figure}

In this case, the bijection \( \overline{\gamma} \) that minimizes the interleaving distance instead matches \( (0,1) \) from the diagram of \( f \) with \( (100,101) \) from the diagram of \( g \), and vice versa. Indeed, this latter leads to a product module isomorphic to \( M_{(f,g)} \).
\end{example}

Our purpose for introducing the $\gamma$-product was to explore how a bigraded module can be constructed from two monoparameter persistence modules. However, it allows us to characterize hook-decomposable persistence modules.

\begin{definition}[\cite{botnan2024signed, botnan2024bottleneck}]
    A biparameter persistence module is called hook-decomposable if it decomposes as a direct sum of hook modules, which are interval modules with a support of the following type:
        \[
             \{ \mathbf{r} \in \mathbb{N}^2 \mid \mathbf{r} \geq \mathbf{p},\ \mathbf{r} \not\geq \mathbf{q} \},
        \]
    for some \( \mathbf{p}, \mathbf{q} \in \mathbb{N}^2 \cup \{\infty\} \) with \( \mathbf{p} \leq \mathbf{q} \) and where  \(\mathbb{N}^2\) is endowed with the structure of product poset.
\end{definition}

Such class of modules is relevant due to the fact that the rank exact resolution of any module always consists of hook-decomposable modules \cite{botnan2024bottleneck}.

By construction, each generator of the \(\gamma\)-product defines a hook summand. Hence, the following result holds.

\begin{proposition}
The module \( M_f \MVAt_{\gamma} M_g \) is hook-decomposable.   
\end{proposition}

Conversely, each hook-decomposable biparameter persistence module coincides with the \(\bar\gamma\)-product of the associated monoparameter persistence modules.

\begin{proposition}
Let \( M_{(f,g)} \) be a hook-decomposable biparameter persistence module. Then, \[ M_{(f,g)} \cong M_f \MVAt_{\bar\gamma} M_g. \]
\end{proposition}




\section*{Acknowledgements}

This work was carried out within the framework of the projects ``RAISE - Robotics and AI for Socio-economic Empowerment'' - Spoke number 2 (Smart Devices and Technologies for Personal and Remote Healthcare), the CNR research activities STRIVE DIT.AD022.207.009 and DIT.AD021.080.001, and the projects “National centre for HPC, Big Data and Quantum Computing” HPC – Spoke 8 and “Molecular imaging study of the immune response in muscle denervation: a high-tech study in murine models and human patients”, and has been supported by European Union - NextGenerationEU.

\bibliographystyle{plain}
\bibliography{main}

\end{document}